%% Amstex file, this is the version of July 16-22, 2016.
%% Corrections July 24-30, 2016, new material  on August 1-2, 2016.
%% August 21, 2016, revision August 28-31, 2016. Major revision
%% September 28, 2018.
%% Minor corrections December 5, 2016. Final version accepted by the "Moscow
%% J. Combinatorics and Number Theory".

\magnification1200
\input amstex.tex
\documentstyle{amsppt}
\nopagenumbers
\hsize=12.5cm
\vsize=18cm
\hoffset=1cm
\voffset=2cm

\footline={\hss{\vbox to 2cm{\vfil\hbox{\rm\folio}}}\hss}

\def\DJ{\leavevmode\setbox0=\hbox{D}\kern0pt\rlap
{\kern.04em\raise.188\ht0\hbox{-}}D}

\def\txt#1{{\textstyle{#1}}}
\baselineskip=13pt
\def\hf{{\textstyle{1\over2}}}
\def\a{\alpha}
\def\d{{\,\roman d}}
\def\e{\varepsilon}\def\E{{\roman e}}

 \def\g{\gamma}
\def\G{\Gamma}

\def\s{\sigma}
\def\t{\theta}
\def\={\;=\;}

\def\zt{\zeta(\hf+it)}

\def\E{{\roman e}}

\def\z{\zeta}

 \def\t{\theta}
\def\hf{{\textstyle{1\over2}}}
\def\txt#1{{\textstyle{#1}}}

\def\le{\leqslant} \def\ge{\geqslant}
%%%%%%%%%%% Fonts macros %%%%%%%%%%%%
\font\tenmsb=msbm10
\font\sevenmsb=msbm7
\font\fivemsb=msbm5
\newfam\msbfam
\textfont\msbfam=\tenmsb
\scriptfont\msbfam=\sevenmsb
\scriptscriptfont\msbfam=\fivemsb
\def\Bbb#1{{\fam\msbfam #1}}

\def \NN {\Bbb N}
\def \CC {\Bbb C}

\font\ff=cmr8
\def\txt#1{{\textstyle{#1}}}
\baselineskip=13pt

\font\teneufm=eufm10
\font\seveneufm=eufm7
\font\fiveeufm=eufm5
\newfam\eufmfam
\textfont\eufmfam=\teneufm
\scriptfont\eufmfam=\seveneufm
\scriptscriptfont\eufmfam=\fiveeufm
\def\mathfrak#1{{\fam\eufmfam\relax#1}}

\font\tenmsb=msbm10
\font\sevenmsb=msbm7
\font\fivemsb=msbm5
\newfam\msbfam
     \textfont\msbfam=\tenmsb
      \scriptfont\msbfam=\sevenmsb
      \scriptscriptfont\msbfam=\fivemsb
\def\Bbb#1{{\fam\msbfam #1}}

\def \NN {\Bbb N}
\def \CC {\Bbb C}

  \def\rightheadline{{\hfil{\ff
  Moments of Hardy's function over short intervals}\hfil\tenrm\folio}}

  \def\leftheadline{{\tenrm\folio\hfil{\ff
   Aleksandar Ivi\'c }\hfil}}
  \def\emptyheadline{\hfil}
  \headline{\ifnum\pageno=1 \emptyheadline\else
  \ifodd\pageno \rightheadline \else \leftheadline\fi\fi}

\topmatter

\title
On certain moments of Hardy's function $Z(t)$ over short intervals
\endtitle
\author   Aleksandar Ivi\'c
 \endauthor

\nopagenumbers

\medskip

\address
Aleksandar Ivi\'c, Katedra Matematike RGF-a
Universiteta u Beogradu, \DJ u\v sina 7, 11000 Beograd, Serbia
\endaddress
\keywords
Riemann zeta-function, Riemann Hypothesis,  Hardy's function, moments, short intervals
\endkeywords
\subjclass
11M06  \endsubjclass

\bigskip
\email {
\tt
aleksandar.ivic\@rgf.bg.ac.rs, aivic\_2000\@yahoo.com }\endemail
\dedicatory
\enddedicatory
\abstract
{Let as usual $Z(t) = \zt\chi^{-1/2}(\hf+it)$ denote Hardy's function,
where $\z(s) = \chi(s)\z(1-s)$. Assuming
the Riemann hypothesis upper and lower bounds for some integrals involving $Z(t)$
and $Z'(t)$ are proved. It is also proved that
$$
H(\log T)^{k^2} \ll_{k,\a} \sum_{T<\g\le T+H}\max_{\g\le \tau_\g\le \g^+}
|\z(\hf + i\tau_\g)|^{2k} \ll_{k,\a} H(\log T)^{k^2}.
$$
Here $k>1$ is a fixed integer, $\g, \g^+$ denote ordinates of consecutive complex zeros of $\z(s)$
 and $T^\a \le H \le T$, where $\a$ is a fixed constant such that $0<\a \le 1$.  This
 sharpens and generalizes a result of M.B. Milinovich \cite{17}.
 }

 \endabstract
\endtopmatter

\document

\head
1. Introduction
\endhead
Let the Riemann zeta-function be, as usual,
$$
\z(s) = \sum_{n=1}^\infty n^{-s}\qquad(\Re s > 1).
$$
For $\Re s \le 1$ one defines $\z(s)$ by analytic continuation (see the monographs
of H.M. Edwards \cite{3}, the author  \cite{11} and E.C. Titchmarsh
 \cite{20} for the properties of $\z(s)$).
Here the Riemann Hypothesis (RH), that all complex
zeros of $\z(s)$ satisfy $\Re s = \hf$, is  assumed {\it throughout the paper}.
The Riemann zeta-function satisfies the functional equation

 $$
\z(s) \;=\; \chi(s)\z(1-s)\quad(\forall s)\in\CC),\quad
\chi(s) := \frac{\G(\hf(1-s))}{\G(\hf s)}\pi^{s-1/2},\leqno(1.1)
$$
where $\G(s)$ is the familiar gamma-function.
 One then defines Hardy's function $Z(t)$ as
 $$
Z(t)  := \zt\bigl(\chi(\hf+it)\bigr)^{-1/2}, \leqno(1.2)
$$
which is real for $t$ real and $|\zt| = |Z(t)|$. Thus the real zeros of $Z(t)$  correspond
to the zeros of $\z(s)$ of the form $\hf+it$, which makes Hardy's function an
invaluable tool in the study of zeros of $\z(s)$ on the critical line $\Re s = \hf$.
For an extensive account on $Z(t)$ the reader is
referred to the author's monograph \cite{16}.

\medskip
Several papers deal with the estimation of the sum
$$
{\Cal M}_k(T): = \frac{1}{N(T)}\sum_{0<\g\le T}\max_{\g\le\tau_\g\le\g^+}|\z(\hf + i\tau_\g)|^{2k}
\equiv \frac{1}{N(T)}\sum_{0<\g\le T}\max_{\g\le\tau_\g\le\g^+}|Z(\tau_\g)|^{2k}.\leqno(1.3)
$$
Here $k\in\NN$ is fixed, and $\g, \g^+$ denote ordinates of consecutive complex zeros of $\z(s)$,
ordered according to their size. Also, as usual,
$$
N(T) = \sum_{0<\g\le T}1 = \frac{T}{2\pi}\log \frac{T}{2\pi} - \frac{T}{2\pi} + O(\log T)\leqno(1.4)
$$
counts (with multiplicities) the number of zeros of $\z(s)$ whose ordinates $\g$ satisfy $0<\g\le T$.

\medskip

In \cite{2} B. Conrey and A. Ghosh proved, under the RH,
$$
{\Cal M}_1(T) = \frac{\E^2-5}{2}\log T + O(1).\leqno(1.5)
$$
Actually, they prove a somewhat stronger result than (1.5), namely
$$
\sum_{T<\g\le T+H}\max_{\g\le\tau_\g\le\g^+}|\z(\hf + i\tau_\g)|^{2}
= \frac{\E^2-5}{4\pi}H\log^2 T + O(H\log T)\leqno(1.6)
$$
with $H = T^{3/4}$. This follows from their proof on noting that (1.4) implies
$$
N(T+H) - N(T) \;\sim\; \frac{H}{2\pi}\log T\qquad(T\to\infty).
$$

B. Conrey \cite{1} obtained, also under the RH,
$$
\frac{\sqrt{21}}{45\pi}(1+o(1))\log^4T \le {\Cal M}_2(T) \le \frac{1+o(1)}{\pi\sqrt{15}}\log^4T
\quad(T\to\infty),\leqno(1.7)
$$
and R.R. Hall \cite{4}, \cite{5} obtained some further improvements of (1.7). A general result,
due to M.B. Milinovich \cite{17}, states that under the RH, for fixed $k\in\NN$.
$$
(\log T)^{k^2-\e} \;\ll_{k,\e}\; {\Cal M}_k(T) \;\ll_{k,\e}\; (\log T)^{k^2+\e}.\leqno(1.8)
$$
Here $\ll_{k,\e}$ means that the  constant implied by the $\ll$-symbol
depends only on $k$ and $\e$, an
arbitrarily small positive number, not necessarily the same one at each occurrence.
The bounds in (1.8), when $k=1,2$, are implied by (1.6) and (1.7), respectively.

\medskip
{\bf Acknowledgment}. The author wishes to thank M.B. Milinovich for valuable remarks.

\head
2. Statement of results
\endhead
\medskip
Milinovich \cite{17} derives (1.8) from upper and lower bounds involving certain
integrals with $Z(t)$ and $Z'(t)$, which seem to be of independent interest.
He investigated integrals over the ``long'' interval $[0, T]$, but here we are
interested in the integrals over the ``short'' intervals $[T, T+H]$, where $H = H(T)$ may
be much smaller than $T$. We shall prove here the following theorems.

\medskip
THEOREM 1. {\it Let $k \ge2$ be a  fixed integer. Under the RH we have, for
$T^\a \le H = H(T) \le T, 0<\a\le1$ a fixed constant},
$$
\int_T^{T+H}(Z'(t))^2Z^{2k-2}(t)\d t \;\gg_{k,\a}\; H(\log T)^{k^2+2}.\leqno(2.1)
$$

\medskip
THEOREM 2. {\it Let $k\in\NN$ be fixed. Under the RH we have, for
$T^\a \le H = H(T) \le T, 0<\a\le1$ a fixed constant,
$$
\int_T^{T+H}{|\z'(\hf+it)|}^{2k}\d t \;\ll_{k,\a}\; H(\log T)^{k^2+2},\leqno(2.2)
$$
and}
$$
\int_T^{T+H}{\bigl(Z'(t)\bigr)}^{2k}\d t \;\ll_{k,\a}\; H(\log T)^{k^2+2}.\leqno(2.3)
$$

\medskip
These bounds differ from the analogous results of \cite{17} in two aspects. Firstly, Milinovich
has the integrals over $[0, T]$, which corresponds to the case $H = T$ in our theorems.
Indeed, if (2.1)--(2.3) hold with $H=T$, then replacing $T$ by $T/2, T/2^2,\ldots$ etc.
and adding up all the results we obtain (2.1)--(2.3) with the interval of integration $[0, T]$.
Secondly,  in (2.1) Milinovich obtained $k^2 +2-\e$ as the exponent of the logarithm,
and in (2.2) and (2.3)
he had the exponents $k^2 +2+\e$ of the logarithm. He remarks on page 1122 that one can
get rid of the $\e$'s in his bounds provided that one has
$$
\int_0^T|\zt|^{2k}\d t \;\ll_k\; T(\log T)^{k^2}\leqno(2.4)
$$
and
$$
\int_0^T|\z'(\hf+it)|^{2k}\d t \;\ll_k\; T(\log T)^{k^2+2k}.\leqno(2.5)
$$
The estimates (2.4) and (2.5) do hold indeed. Namely K. Soundararajan \cite{19} proved
(2.4) with the exponent of $\log T$ in (2.4) equal to $k^2+\e$. The author \cite{15} improved
and sharpened Soundararajan's bound by showing that
$$
 \int\limits_{T}^{T+H}|\zt|^{2k}\d t
\;\ll_{k,\a}\; H(\log T)^{k^2\bigl(1+O(1/\log_3T)\bigr)}\quad(\roman{RH}).\leqno(2.6)
$$
Here $T^\a \le H \le T$ where $0 < \a \leqslant 1$ is a fixed number, and
$$
\log_3T = \log\log\log T = \log(\log_2T).
$$
Note that \cite{15} appeared before \cite{19} because of the
backlog of ``Ann. Math.'' The key result in \cite{15}, which is a proper generalization of
the corresponding result in \cite{19}, is

\medskip
\medskip
THEOREM A. {\it Let $H = T^\t$ where $0 < \t \leqslant 1$ is a fixed number, and
let $\mu(T,H,V)$ denote the measure of points $t$ from $[T-H,\,T+H]$ such that
$$
\log|\zt| \geqslant V,\quad 10\sqrt{\log_2T} \leqslant V \leqslant {3\log 2T\over 8\log_2(2T)}.
%\leqno(7.6)
$$
Then, under the RH, for $\;10\sqrt{\log_2T} \leqslant V \leqslant\log_2T$  we have
$$
\mu(T,H,V) \;\ll\; H{V\over\sqrt{\log_2T}}\exp\left(-{V^2\over\log_2T}\Bigl(1-
{7\over2\t\log_3T}\Bigr)\right),
%\leqno(7.7)
$$
for $\;\log_2T \leqslant V \leqslant\hf\t\log_2T\log_3T$ we have
$$
\mu(T,H,V) \;\ll\; H\exp\left(- {V^2\over\log_2T}
\Bigl(1 - {7V\over4\t\log_2T\log_3T}\Bigr)^2\,\right),
%\leqno(7.8)
$$
and for $\;\hf\t\log_2T\log_3T\leqslant V\leqslant {3\log 2T\over 8\log_2(2T)}$ we have}
$$
\mu(T,H,V) \;\ll\; H\exp(-\txt{1\over20}\t V\log V).
%\leqno(7.9)
$$

Later A. Harper \cite{9} (RH) improved the upper bound in \cite{19}
by establishing  (2.4) for the long interval $[0, T]$. As remarked in \cite{9} on p. 4,
the method of \cite{15} leading to (2.6), i.e., Theorem A,
can be combined with that of \cite{9} to produce the sharp upper
bound over the short interval $[T, T+H]$, namely
$$
 \int\limits_{T}^{T+H}|\zt|^{2k}\d t
\;\ll_{k,\a}\; H(\log T)^{k^2}\quad({\roman{RH})}.\leqno(2.7)
$$
The bound in (2.7) is the key ingredient in the proof of our results.
It is, up to the values of the $\ll$-constants, best possible, since long ago it was shown by
R. Balasubramanian and K. Ramachandra (see the latter's monograph \cite{18}, in particular
the remark on p. 45)
that,
if $k\geqslant1$ is a fixed integer, then for $C(\e,k)\log\log T
\leqslant H \leqslant T/2$ we have
$$
\int\limits_{T}^{T+H}|\zt|^{2k}\d t \;\geqslant\; (C_k'-\e)H(\log H)^{k^2},\leqno(2.8)
$$
where
$$
C_k' = {1\over2\G(k^2+1)}\prod_p\left\{(1-p^{-1})^{k^2}\sum_{m=0}^\infty
\Bigl({\G(k+m)\over\G(k)m!}\Bigr)^2p^{-m}\right\}.
$$
We note that the  lower bound in (2.8) is {\it unconditional}, with a very wide range for $H$.
As for (2.5), it will be shown later that a corresponding result holds over $[T, T+H]$.

\medskip
THEOREM 3. {\it
Let $1< k\in\NN$ be fixed, $\g, \g^+$ denote ordinates of consecutive complex zeros of $\z(s)$
 and $T^\a \le H = H(T)\le T$, where $\a$ is a fixed constant such that $0<\a \le 1$.
 Under the RH we have then
$$
H(\log T)^{k^2} \ll_{k,\a} \sum_{T<\g\le T+H}\max_{\g\le \tau_\g\le \g^+}
|\z(\hf + i\tau_\g)|^{2k} \ll_{k,\a} H(\log T)^{k^2}.\leqno(2.9)
$$
}

\medskip
{\bf Remark}. The case $k=1$ was treated in \cite{2} (see (1.6))
 and is not covered by Theorem 3. This is because Theorem 1
does not cover the case $k=1$. The method of its proof (see (3.4))
does not work in obtaining a lower bound for
$$
\int_T^{T+H}|Z'(t)Z(t)|\d t.
$$

\medskip
\head
3. Proof of Theorem 1
\endhead
\medskip
As was also done in \cite{Mil}, we follow Conrey and Ghosh \cite{2}, and introduce the
analytic function
$$
Z_1(s) := \z'(s) - \frac{\chi'(s)}{2\chi(s)}\z(s).\leqno(3.1)
$$
Its usefulness  comes from the fact that differentiation of (1.2) gives
$$
Z'(t) = i\left\{\z'(\hf+it) - \hf \frac{\chi'(\hf+it)}{\chi(\hf+it)}
\zt\right\}\chi^{-1/2}(\hf+it).
\leqno(3.2)
$$
This gives
$$
|Z'(t)| \= |Z_1(\hf+it)|.\leqno(3.3)
$$
Using (3.3) and $|Z(t)| = |\zt|$ one may write
$$
\int_T^{T+H}(Z'(t))^2Z^{2k-2}(t)\d t =
\int_T^{T+H}\Bigl|Z_1(\hf+it)\z(\hf+it)^{k-1}\Bigr|^2\d t.\leqno(3.4)
$$
The basic idea is to use the inequality
$$
\eqalign{&
\Bigl|\int_T^{T+H}Z_1(\hf+it)\z(\hf+it)^{k-1}{\bar A}(t)\d t\Bigr|^2
\cr&
\le \int_T^{T+H}(Z'(t))^2Z^{2k-2}(t)\d t\cdot \int_T^{T+H}|A(t)|^2\d t.\cr}\leqno(3.5)
$$
This comes on using (3.4) and the Cauchy-Schwarz inequality for integrals with a suitably
chosen function $A(t)$. Following \cite{Mil} we set
$$
A(t):= {\Cal A}(\hf + it), \quad {\Cal A}(s) = {\Cal A}(s;k,\xi) := \sum_{n\le\xi}d_k(n)n^{-s},
$$
where $d_k(n)$ (generated by $\z^k(s)$ for $\Re s > 1)$ is the (generalized) divisor
function which represents the number of ways $n$ can be written as a product of $k$ fixed factors
(see e.g., Chapter 13 of \cite{11} for more properties). The parameter $\xi$ is given by
$\xi = T^\t, 0<\t < 1$. The function $A(t)$ has the property that in mean square it behaves like
$(\log \xi)^{k^2}$ (see Chapter 13 of \cite{11}). Therefore, by the mean value theorem
for Dirichlet polynomials (see e.g., Theorem 5.2 of \cite{11}), we have
$$
\eqalign{
\int_T^{T+H}|A(t)|^2\d t &= H\sum_{n\le\xi}d_k^2(n)n^{-1} + O\Bigl(\sum_{n\le\xi}d_k^2(n)\Bigr)
\cr&
= H(C_k + o(1))(\log\xi)^{k^2} + O(\xi(\log\xi)^{k^2-1})\cr}\leqno(3.6)
$$
for $2\le \xi \le T$ and a positive constant $C_k$, which may be made explicit.
It remains to estimate from below
$$
\eqalign{&
\int_T^{T+H}Z_1(\hf+it)\z(\hf+it)^{k-1}{\bar A}(t)\d t
\cr&
= \frac{1}{i}\int_{1/2+iT}^{1/2+iT+iH}Z_1(s)\z^k(s){\Cal A}(1-s)\d s
\cr&
= \int_{a+iT}^{a+iT+iH}\frac{1}{i}Z_1(s)\z^k(s){\Cal A}(1-s)\d s + O_\e(T^\e\xi).\cr}
\leqno(3.7)
$$
Here we used Cauchy's theorem and set $a := 1 + 1/\log T$. We also used
standard consequences of the RH  (see Chapter 12 of \cite{20}):
$$
\z(s) \ll_{\e,\s} |t|^\e, \quad\z'(s) \ll_{\e,\s} |t|^\e\qquad(s = \s+it, \s \ge\hf),
$$
as well as the unconditional, elementary bound $d_k(n) \ll_{\e,k} n^\e$ and
$$
\frac{\chi'(\s+it)}{\chi(\s+it)} = -\log\frac{t}{2\pi} + O\Bigl(\frac{1}{t}\Bigr).
\leqno(3.8)
$$
One obtains (3.8) by logarithmic differentiation of (1.1) and the use of Stirling's
formula for the gamma-function. It is valid for $\hf\le\Re s \le 2$, and the $O$-term
in (3.8) in fact admits a full asymptotic expansion in terms of negative exponents
of $t$, and the left-hand side of (3.8) can be further differentiated.
The integral on the right-hand side of (3.7) is written as $J_1+J_2$, where
$$
\eqalign{&
J_1: \= \frac{1}{i}\int_{a+iT}^{a+iT+iH}\z'(s)\z^{k-1}{\Cal A}(1-s)\d s,\cr&
J_2: \= -\frac{1}{2i}\int_{a+iT}^{a+iT+iH}\frac{\chi'(s)}{\chi(s)}\z^k(s){\Cal A}(1-s)\d s.
\cr}
$$
Similarly as in \cite{17} one shows that
$$
J_1 = -H\sum_{n\le\xi}{\tilde d}_k(n)d_k(n)n^{-1} + O(T^\e\xi)
$$
with
$$
{\tilde d}_k(n): = \sum_{\delta|n}d_{k-1}(\delta)\log\frac{n}{\delta} \le d_k(n)\log n,
$$
and
$$
J_2 = \hf H\left(\log\frac{T}{2\pi{\roman e}}\right)\sum_{n\le\xi}d_k^2(n)n^{-1}+ O(T^\e\xi).
$$
This yields
$$
\eqalign{&
\int_T^{T+H}Z_1(\hf+it)\z(\hf+it)^{k-1}{\bar A}(t)\d t = J_1+J_2 + O(T^\e\xi)\cr&
= -H\sum_{n\le\xi}{\tilde d}_k(n)d_k(n)n^{-1} + O(T^\e\xi)
+ \hf H\log {\frac{T}{2\pi}}(1+o(1))\sum_{n\le\xi}d_k^2(n)n^{-1}\cr&
\ge -H\sum_{n\le\xi}d_k^2(n)n^{-1}\log n +
\hf H\log {\frac{T}{2\pi}}(1+o(1))\sum_{n\le\xi}d_k^2(n)n^{-1}+  O(T^\e\xi)\cr&
\ge \Bigl\{\hf H(1+o(1))\log\frac{T}{2\pi} - H\log\xi\Bigr\}\sum_{n\le\xi}d_k^2(n)n^{-1}
+ O(T^\e\xi)\cr&
\ge A_k H\log T\cdot(\log\xi)^{k^2}\cr}
$$
for $\xi = T^{\t}, \t = \hf \a$ and $\e$ sufficiently small. From (3.5) and (3.6) we finally gather that
$$
H^2\log^2T(\log T)^{2k^2} \ll_{k,\a} \int_T^{T+H}(Z'(t))^2Z^{2k-2}(t)\d t\cdot H(\log T)^{k^2},
$$
and (2.1) of Theorem 1 follows.

\medskip
\head
4. Proof of Theorem 2
\endhead

\medskip

First note that, for $0<R\le\hf, T\le t \le 2T$, by Cauchy's integral formula we have
$$
\z'(\hf+it) \;=\;\frac{1}{2\pi i}\int_{|z|=R}\frac{\z(\hf+it+z)}{z^2}\d z.
$$
This yields
$$
\int_T^{T+H}|\z'(\hf+it)|^{2k}\d t = \frac{1}{(2\pi)^{2k}}\int_T^{T+H}
\left|\int_{|z|=R}\frac{\z(\hf+it+z)}{z^2}\d z\right|^{2k}\d t.
\leqno(4.1)
$$
By H\"older's inequality for integrals the right-hand side of (4.1) does not exceed
$$
\eqalign{&
\frac{1}{(2\pi)^{2k}}\int_T^{T+H}\left\{\int_{|z|=R}|\z(\hf+it+z)|^{2k}|\d z|\right\}
\cdot
\left\{\frac{|\d z|}{|z|^{4k/(2k-1)}}\right\}^{2k-1}\d t
\cr&
\le \frac{1}{(2\pi)^{2k}}\int_T^{T+H}\left\{\int_{|z|=R}|\z(\hf+it+z)|^{2k}|\d z|\right\}
(2\pi R)^{2k-1}R^{-4k}
\cr&
\le \frac{1}{ R^{2k}}\max_{0\le\t\le2\pi}\int_T^{T+H}|\z(\hf+it + R\E^{i\t})|^{2k}\d t.
\cr}
$$
Therefore
$$
\int_T^{T+H}|\z'(\hf+it)|^{2k}\d t \le
\frac{1}{ R^{2k}}\max_{0\le\t\le2\pi}\int_T^{T+H}|\z(\hf+it + R\E^{i\t})|^{2k}\d t.
\leqno(4.2)
$$
As in \cite{Mil}, we could have obtained an inequality for the $\ell$-th derivative
of $\zt$, but this is not necessary for our purposes.

\medskip
Henceforth let $R=1/\log T$ in (4.2). The integral on the right-hand side of (4.2) equals
$$
\int_T^{T+H}|\z\bigl(\hf+R\cos\a + i(t + R\sin\a)\bigr)|^{2k}\d t.\leqno(4.3)
$$
Recall that, under the RH (see \cite {Tit}),
$$
\z(\s+it) \;\ll\; \exp\left(C\frac{\log t}{\log\log t}\right)\qquad(\hf \le \s\le1, C>0, |t|\ge 2).\leqno(4.4)
$$
When $\cos\a\ge0$ in (4.3), we use (4.4) to obtain that the integral in (4.3) is equal to
$$
\int_T^{T+H}|\z(\hf+R\cos\a + it)|^{2k}\d t + o(T).\leqno(4.5)
$$
When $\pi/2 \le \t \le 3\pi/2$ in (4.5), that is, when $\cos\a \le 0$, we use
the functional equation (1.1). In this case we have, with $\s=\Re s = \hf + R\cos\a, R = 1/\log T$,
$$
\chi(s) \ll |t|^{1/2-\s} \ll T^{-\cos\a/\log T} \ll 1,
$$
thus we reduce the estimation of our integral  to the case when $\cos\a\ge0$.
For this we shall use a convexity result which shows that essentially the integral
in question is bounded by the $2k$-th moment of $|\zt|$ over a short interval.
More precisely, let $1/2\le \s \le 3/4, k>0, t\ge2$,
$$
J_k(\s) :=\int_{-\infty}^\infty|\z(\s+it)|^{2k}w_k(t)\d t,
\quad w_k(t) : =\int_T^{T+H}\E^{-2k(t-\tau)^2}\d\tau.\leqno(4.6)
$$
Then
$$
J_k(\s) \ll T^{\s-1/2}\Bigl(J_k(\hf)\Bigr)^{3/2-\s} + \E^{-kT^2/4}.\leqno(4.7)
$$
The bound in (4.7) is the analogue of Lemma 4.2 of \cite{17} for short intervals. This in turn is a
result of D.R. Heath-Brown \cite{10}, which is also expounded in \cite{12}, pp. 321-323.
In the original version the interval of integration in the kernel function function $w_k(t)$
was $[T, 2T]$. However, the change made in (4.6) does not affect the proof, and one obtains (4.7).
Now note that $w_k(t) \gg 1$ for $t\in [T,T+H]$, so that
$$
\int_T^{T+H}|\z(\s+it)|^{2k}\d t \ll J_k(\s).\leqno(4.8)
$$
We have the bound $w_k(t) \ll \exp(-2kH^2)$ when $t \le T-H$ or $t \ge T+2H$. On the other hand
$w_k(t) \ll \exp(-kt^2)$ for $t<0$ or $t > 3T$. Thus combining (4.7) and (4.8) it follows that
$$
\eqalign{
\int_T^{T+H}|\z(\s+it)|^{2k}\d t &= \int_{-\infty}^0 + \int_0^{T-H} + \int_{T-H}^{T+2H} + \int_{T+2H}^\infty
\cr&
\ll 1 + \int_{T-H}^{T+2H}|\z(\hf+it)|^{2k}\d t\cr& \ll_{k,\a} H(\log T)^{k^2},
\cr}\leqno(4.9)
$$
where in the last step (2.7) was used.
Inserting (4.9) in (4.2) the bound in (2.2) follows. The estimate (2.3) easily
follows from (2.2), (2.9) and (3.2).
Theorem 2 is proved.

\medskip
\head
5. Proof of Theorem 3
\endhead

\medskip
It is in the folklore that $Z(t)$, for $t\ge14$, cannot have a negative local
maximum or a positive local minimum under the RH. For this, see \cite{3}, or
\cite{13}, \cite{14}, \cite{6}. In other words, the zeros of $Z(t)$ and $Z'(t)$ are
interlacing. Thus if $\g, \g^+$ are consecutive zeros of $Z(t)$, there is a unique point
$\lambda_\g\in [\g, \g^+]$ for which $Z'(\lambda_\g)=0$ (this is trivially true if
$\g=\g^+$, that is, if $\g$ is a multiple zero of $Z(t)$). Therefore
$$
\max_{\g\le \tau_\g\le \g^+}|\z(\hf + i\tau_\g)|^{2k} = Z^{2k}(\lambda_\g).\leqno(5.1)
$$
Then, since $Z(t)$ is positive in $(\g, \lambda_\g)$ and negative in $(\lambda_\g, \g^+)$,
or conversely, we have
$$
\eqalign{
\int_\g^{\g^+}|Z'(t)Z^{2k-1}(t)|\d t &=
\left|\int_\g^{\lambda_\g}Z'(t)Z^{2k-1}(t)\d t - \int_{\lambda_\g}^\g Z'(t)Z^{2k-1}(t) \d t
\right|
\cr&
= \frac{1}{k}Z^{2k}(\lambda_\g).\cr}\leqno(5.2)
$$
Therefore (5.1) and (5.2) give, in view of (4.4),
$$
\sum_{T<\g\le T+H}\max_{\g\le \tau_\g\le \g^+}|\z(\hf + i\tau_\g)|^{2k}
= k\int_T^{T+H}|Z'(t)Z^{2k-1}(t)|\d t + O_{k,\e}(T^\e).\leqno(5.3)
$$

Assume  that $k\ge2$, so that $2k-2\ge2$. To bound
the integral on the right-hand side of  (5.3) from below, note that
$$
|(Z')^2Z^{2k-2}| = |Z'|^{1/2}|Z|^{k-1/2}\cdot |Z'|^{3/2}\cdot|Z|^{k-3/2}.
$$
Thus H\"older's inequality for integrals shows that
$$
\eqalign{&
\int_T^{T+H}(Z'(t))^2Z^{2k-2}(t)\d t\le
\cr&
 \left(\int_T^{T+H}(|Z'|^{1/2}|Z|^{k-1/2})^p\d t\right)^{\frac{1}{p}}
\left(\int_T^{T+H}|Z'|^{3q/2}\d t\right)^{\frac{1}{q}}
\left(\int_T^{T+H}|Z|^{r(k-3/2)}\d t\right)^{\frac{1}{r}}\cr}
$$
with $p,q,r>0, 1/p+1/q+1/r=1$. Take
$$
\frac{1}{p} = \frac{1}{2},\quad \frac{1}{q} = \frac{3}{4k},
\quad \frac{1}{r} = \frac{1}{2} - \frac{3}{4k}.
$$
Then the  right-hand side is, on using (2.3) and (2.4),
$$
\eqalign{&
\le \left(\int_T^{T+H}|ZZ|^{2k-1}\d t\right)^{\frac{1}{2}}\left(\int_T^{T+H}|Z'|^{2k}\d t\right)^{\frac{3}{4k}}
\left(\int_T^{T+H}|Z|^{2k}\d t\right)^{\frac{1}{2}-\frac{3}{4k}}
\cr&
\ll_{k,\a}  \left(\int_T^{T+H}|ZZ|^{2k-1}\d t\right)^{\frac{1}{2}}\left(H(\log T)^{k^2+2k}\right)^{\frac{3}{4k}}
\left(H(\log T)^{k^2}\right)^{\frac{1}{2}-\frac{3}{4k}}.\cr}
$$
This gives, on using (2.1),
$$
H(\log T)^{k^2+2} \;\ll_{k,\a}\;
 I^{1/2} \left(H(\log T)^{k^2+2k}\right)^{\frac{3}{4k}}\left(H(\log T)^{k^2}\right)^{\frac{1}{2}-\frac{3}{4k}},
$$
which on simplifying yields
$$
I :=  \int_T^{T+H}|Z'(t)Z^{2k-1}(t)|\d t  \;\gg_{k,\a}\; H(\log T)^{k^2+1}.
$$
In view of (5.3) this proves the lower bound in (2.9) of Theorem 3.

 As for the upper bound, the integral in (5.3)
does not exceed, by H\"older's inequality for integrals,
$$
\eqalign{&
\left|\int_T^{T+H}(Z'(t))^{2k}\d t\right|^{1/(2k)}
\left|\int_T^{T+H}|\zt|^{2k}\d t\right|^{1-1/(2k)}
\cr&
\ll_{k,\a} \left\{H(\log T)^{k^2+2k}\right\}^{1/(2k)}
\left\{H(\log T)^{k^2}\right\}^{1-1/(2k)}
\cr&
= H(\log T)^{k^2+1},\cr}\leqno(5.4)
$$
which finishes the proof of Theorem 3 when $k\ge2$. Here we used (2.3) and (2.4), and we note that
the bound in (5.4) holds also for $k=1$. It is the lower bound in this case which
is problematic.

%\vfill
%\eject
%\topglue1cm
\bigskip
\Refs
\bigskip

\item{[1]} J.B. Conrey, The fourth moment of derivatives of the Riemann
zeta-function, Quart. J. Math. Oxford Ser. (2){\bf39}(1988), 21-36.

\item{[2]} J.B. Conrey and A. Ghosh, A mean value theorem for the Riemann zeta-function
at its relative extrema on the critical line,
 J. Lond. Math. Soc., II. Ser. {\bf32}(1985), 193-202.

\item{[3]} H.M. Edwards, Riemann's zeta-function, Academic Press, New York-London, 1974.

\item{[4]} R.R. Hall, The behaviour of the Riemann zeta-function
on the critical line, Mathematica {\bf46}(1999), 281-313.

\item{[5]} R.R. Hall, On the extreme values of the Riemann zeta-function
between its zeros on the critical line, J. reine angew. Math. {\bf560}(2003),
29-41.

\item{[6]} R.R. Hall, On the stationary points of Hardy's function $Z(t)$,
Acta Arith. {\bf111}(2004), 125-140.

\item{[7]} R.R. Hall, A new unconditional result about large spaces
between zeta zeros, Mathematika {\bf53}(2005), 101-113.

\item{[8]} R.R. Hall, Extreme values of the Riemann zeta-function
on short zero intervals, Acta Arith. {\bf121}(2006), 259-273.

\item{[9]} A.J. Harper, Sharp conditional bounds for moments of the Riemann zeta-function,
preprint available at {\tt arXiv:1305.4618}.

\item{[10]} D.R. Heath-Brown, Fractional moments of the Riemann zeta-function, J. London
Math. Soc. {\bf24}(1981), 65-78.

\item{[11]} A. Ivi\'c, The Riemann zeta-function, John Wiley \&
Sons, New York 1985 (reissue,  Dover, Mineola, New York, 2003).

\item {[12]} A. Ivi\'c,  Mean values of the Riemann zeta-function,
LN's {\bf 82},  Tata Inst. of Fundamental Research,
Bombay,  1991 (distr. by Springer Verlag, Berlin etc.).

\item {[13]} A. Ivi\'c, On some results concerning the Riemann Hypothesis, in
``Analytic Number Theory" (Kyoto,
1996) ed. Y. Motohashi, LMS LNS {\bf247}, Cambridge University Press,
Cambridge, 1997, pp. 139-167.

\item {[14]} A. Ivi\'c,
On some reasons for doubting the Riemann Hypothesis,
in P. Borwein, S. Choi, B. Rooney and A. Weirathmueller, ``The Riemann
Hypothesis'', CMS Books in Mathematics, Springer, 2008.

\item {[15]} A. Ivi\'c, On mean value results
for the Riemann Zeta-Function in short intervals,
Hardy-Ramanujan J. {\bf32}(2009), 4-23.

\item{[16]} A. Ivi\'c,  The theory of Hardy's $Z$-function,
Cambridge University Press, Cambridge, 2012, 245pp.

\item{[17]} M.B. Milinovich, Moments of the Riemann zeta-function at its
relative extrema on the critical line, Bull. London Math. Soc. {\bf43}(2011), 1119-1129.

\item{[18]} K.  Ramachandra, On the mean-value and omega-theorems
for the Riemann zeta-function, LN's {\bf85}, Tata Inst. of Fundamental Research
(distr. by Springer Verlag, Berlin etc.), Bombay, 1995.

\item{[19]} K. Soundararajan, Moments of the Riemann zeta function, Ann. Math.
(2) {\bf170}(2009), 981-993.

\item{[20]}  E.C. Titchmarsh, The theory of the Riemann
zeta-function (2nd edition),  Oxford University Press, Oxford, 1986.

\endRefs

\enddocument

\bye